\documentclass{article}%
\usepackage{amsmath}
\usepackage{amsfonts}
\usepackage{amssymb}
\usepackage{graphicx}%
\providecommand{\U}[1]{\protect \rule{.1in}{.1in}}
\newtheorem{theorem}{Theorem}

\newtheorem{corollary}[theorem]{Corollary}

\newtheorem{definition}[theorem]{Definition}

\newtheorem{lemma}[theorem]{Lemma}

\newtheorem{proposition}[theorem]{Proposition}
\newtheorem{remark}[theorem]{Remark}

\newenvironment{proof}[1][Proof]{\noindent \textbf{#1.} }{\  \rule{0.5em}{0.5em}}

\def\essinf{\operatorname{ess.\!inf}}
\def\esssup{\operatorname{ess.\!sup}}
\def\sintr{]\!]}
\def\sintl{[\![}

\begin{document}

\title{\bf Representation of the penalty term of dynamic concave utilities}

\author{Freddy Delbaen \footnote{Department of Mathematics, ETH, Z\"urich, Switzerland.
e-mail: delbaen@math.ethz.ch} $^,$ \thanks{This research was
sponsored by a grant of Credit Suisse as well as by a grant
NCCR-Finrisk. The text only reflects the opinion of the authors.
\qquad \qquad \break Part of the research was done while this author
was visiting China in 2005, 2006 and 2007. The hospitality of
Shandong University is greatly appreciated.}, \quad Shige Peng
\footnote{Institute of Mathematics, Shandong University, Jinan,
China. e-mail: peng@sdu.edu.cn} \quad and \quad Emanuela Rosazza
Gianin \footnote{Corresponding author.
Dipartimento di Metodi Quantitativi
per le Scienze Economiche ed Aziendali, Universit\`{a} di
Milano-Bicocca, Italy. e-mail: emanuela.rosazza1@unimib.it} $^,$
\thanks{This work was done while this author was appointed at the Universit\`{a} di Napoli ``Federico II'', Italy. Part of this research was carried out during her visiting in
China in 2006 and at ETH in Z\"urich in 2004, 2006 and 2007. The warm
hospitality of Shandong University and of ETH is gratefully
acknowledged.} }

\date{November 4, 2008}

\maketitle

\abstract{In the context of a Brownian filtration and with a fixed
finite time horizon, we will provide a representation of the penalty
term of general dynamic concave utilities (hence of dynamic convex
risk measures) by applying the theory of $g$-expectations.

\section{Introduction}

Coherent risk measures were introduced by Artzner et al. \cite{ADEH}
in finite sample spaces and later by Delbaen \cite{Delb1} and
\cite{Delb2} in general probability spaces. The aim of this
financial tool is to quantify the intertemporal riskiness which an
investor would face at a maturity date $T$ in order to decide if
this risk could be acceptable for him or not. The family of coherent
risk measures were extended later by F{\" o}llmer and Schied
\cite{FS1}, \cite{FS2} and Frittelli and Rosazza Gianin
\cite{FrRos1}, \cite{FrRos2} to the class of convex risk measures.

$g-$expectations were introduced by Peng \cite{Peng1} as solutions
of a class of nonlinear Backward Stochastic Differential Equations
(BSDE, for short), a class which was first studied by Pardoux and
Peng \cite{PP}. Financial applications and particular cases were
discussed in detail by El Karoui et al. \cite{EPQ}.

As shown by Rosazza Gianin \cite{Ros}, the families of static risk
measures and of $g-$expectations are not disjoint. Indeed, under
suitable hypothesis on the functional $g$, $g-$expectations provide
examples of coherent and/or convex static risk measures.
Furthermore, by defining ``dynamic risk measure" as a ``map" which
quantifies at any intermediate time $t$ the riskiness which will be
faced at maturity $T$, a class of dynamic risk measures can be
obtained by means of conditional $g-$expectations. In particular,
any dynamic risk measure induced by a conditional $g-$expectation
satisfies a ``time-consistency property" (in line with the notion
introduced by Koopmans \cite{Koop} and Duffie and Epstein
\cite{DuffieEp}) or, in the language of Artzner et al.
\cite{ADEHK2}, a ``recursivity property". Further discussions on
dynamic risk measures can be found in Artzner et al. \cite{ADEHK2},
Barrieu and El Karoui \cite{BaElK2}, Bion-Nadal \cite{BionN},
Cheridito et al. \cite{CDK1}, \cite{CDK2}, Cheridito and Kupper
\cite{CK}, Detlefsen and Scandolo \cite{DetlS}, Frittelli and
Rosazza Gianin \cite{FrRos2} and Kl\"{o}ppel and Schweizer
\cite{KlS}, among many others.

\medskip
The main aim of this paper is to represent the penalty term of
general dynamic concave utilities (hence of dynamic convex risk
measures) in the context of a Brownian filtration, a fixed finite
time horizon $T$ and under the assumption of the existence of an
equivalent probability measure with zero penalty. By applying the
theory of $g$-expectations, we will finally prove that the penalty
term is of the following form:
$$
c_{s,t}(Q)=E_Q \left[\int _s ^t f(u,q_u) du |
\mathcal{F}_s
 \right]
$$
(see the exact statement in Theorem \ref{thm3}).

\medskip
The paper is organised as follows. Some well-known results on BSDE
and on risk measures are recalled in Section \ref{section2}. Section
\ref{section3} contains the main result of the paper, that is the
representation of the penalty term of suitable dynamic concave
utilities. As we will see later, this representation will be
obtained by applying the theory of $g$-expectations.

\section{Notation and preliminaries} \label{section2}

Let $(B_{t})_{t\geq0}$ be a standard $d-$dimensional Brownian motion defined
on the probability space $(\Omega,\mathcal{F},P)$ and $\{ \mathcal{F}%
_{t}\}_{t\geq0}$ be the augmented filtration generated by $(B_{t})_{t\geq0}$.

In the sequel, we will identify a probability measure $Q \ll P$ with
its Radon-Nykodim density $\frac{dQ}{dP}$. Furthermore, because of
the choice of the Brownian setting, we will also identify a
probability measure $Q$ equivalent to $P$ with the predictable
process $(q_t)_{t \in [0,T]}$ induced by the stochastic exponential,
i.e. such that
\begin{equation} \label{eq 11}
E_{P}\left[  \frac{dQ}{dP}| \mathcal{F}_{t}\right]
=\mathcal{E} (q.B)_{t}\triangleq \operatorname{exp}\left(
-\frac{1}{2}\int_{0}^{t}\Vert
q_{s}\Vert^{2}ds+\int_{0}^{t}q_{s}dB_{s}\right)
\end{equation}
(see Proposition VIII.1.6 of Revuz and Yor \cite{RevuzYor}).

\bigskip
Consider now a function
$$
\begin{array}{rlclclclll}
g:&{\Bbb R}^{+} &\times &\Omega &\times &\Bbb R &\times &{\Bbb R}^{d} &\rightarrow &\Bbb R\\
&(t, & &\omega ,& &y, & &z) &\longmapsto &g(t,\omega ,y,z)%
\end{array}
$$
satisfying at least the following assumptions (as in Coquet et al.
\cite{CHMP}, but without imposing a priori an horizon of time $T$).
To simplify the notations, we will often write $g(t,y,z)$ instead of
$g(t,\omega,y,z)$.
\bigskip

\textbf{Basic assumptions on $g$:}

(A) $g$ is Lipschitz in $(y,z)$, i.e. there exists a constant $\mu>0$ such
that, $(dt\times dP)-a.s.$, for any $(y_{0},z_{0}),(y_{1},z_{1})\in
\mathbb{R}\times \mathbb{R}^{d}$,
\[
\vert g(t,y_{0},z_{0})-g(t,y_{1},z_{1})\vert \leq \mu(\vert y_{0}-y_{1}%
\vert+\Vert z_{0}-z_{1}\Vert).
\]

(B) For all $(y,z) \in \Bbb R \times \Bbb R ^d$, $g(\cdot,y,z)$ is a
predictable process such that for any finite $T>0$ it holds
$E[\int_{0}^{T}(g(t,\omega,y,z))^{2}dt]<+\infty$ for any $y\in
\mathbb{R}$ and $z\in \mathbb{R}^{d}$.
\smallskip

(C) $(dt\times dP)-$a.s., $\forall y\in \mathbb{R}$, $g(t,y,0)=0$.
\bigskip

Once the horizon of time $T>0$ is fixed, Pardoux and Peng \cite{PP}
introduced the following Backward Stochastic Differential Equation
(BSDE, for short):
$$
\left\{
\begin{array}{rl}
-dY_{t}&=g(t,Y_{t},Z_{t})dt-Z_{t}dB_{t}\\
Y_{T}&=\xi,
\end {array} \right.
$$
where $\xi$ is a random variable in
$L^{2}(\Omega,\mathcal{F}_{T},P)$. Moreover, they showed (see also
El Karoui et al. \cite{EPQ}) that there exists a unique solution
$(Y_{t},Z_{t})_{t\in[0,T]}$ of predictable stochastic processes (the
former $\mathbb{R}$-valued, the latter $\mathbb{R}^{d}$-valued) such
that $E[\int_{0}^{T}Y_{t}^{2}dt]<+\infty$ and $E[\int_{0}^{T}\Vert
Z_{t}\Vert^{2}dt]<+\infty$.

Peng \cite{Peng1} defined the $g-$expectation of $\xi$ as:
$$\mathcal{E}_{g}(\xi)\triangleq Y_{0}$$
and the conditional $g-$expectation of $\xi$ at time $t$ as:
$$\mathcal{E}_{g}(\xi \vert \mathcal{F}_{t})\triangleq Y_{t}.$$
When $g(t,y,z)=\mu \Vert z \Vert$ (with $\mu >0$), $\mathcal{E}_g$
will be denoted by $\mathcal{E}^{\mu}$.
\smallskip

In the sequel, we will only consider essentially bounded random
variables $\xi$, i.e. $\xi \in L^{\infty} (\Omega,\mathcal{F}_T,P)$.

\bigskip

\textbf{Further assumptions on $g$}\smallskip

(1$_{g}$) $g$ does not depend on $y$ \smallskip

(2$_{g}$) $g$ is convex in $z$:

$\qquad \forall \alpha \in [0,1],\forall z_{0},z_{1}\in \mathbb{R}^{d}%
,\quad(dt\times dP)-$a.s.:

$\qquad g(t,\alpha z_{0}+(1-\alpha)z_{1})\leq \alpha
g(t,z_{0})+(1-\alpha)g(t,z_{1})$. \medskip

In the sequel, we will write  {\it \textquotedblleft$g$ with the
usual assumptions"} when $g$ satisfies hypothesis (A)-(C) and
(1$_{g}$)-(2$_{g}$).

\bigskip
Some sufficient conditions for a functional to be induced by a
$g$-expectation are provided by Coquet et al. \cite{CHMP}. Before
recalling this result, we will introduce what is needed.

\begin{definition} \textbf{(Coquet et al.
\cite{CHMP})} A functional $\mathcal{E}: L^2 (\mathcal{F}_T) \to
\Bbb R$ is called an \emph{$\mathcal{F}$-consistent expectation} if
it satisfies the following properties:
\begin{itemize}
\item[(i)] \emph{constancy}: $\mathcal{E}(c)=c$, for any $c \in \Bbb R$;

\item[(ii)] \emph{strict monotonicity}: if $\xi \geq \eta$, then
$\mathcal{E}(\xi) \geq \cal{E}(\eta)$. Moreover, if $\xi \geq \eta$:
$\xi=\eta$ if and only if $\mathcal{E}(\xi)=\mathcal{E}(\eta)$;

\item[(iii)] \emph{consistency}: for any $\xi \in L^2 (\mathcal{F}_T)$ and $t \in [0,T]$
there exists a random variable $\mathcal{E}(\xi \vert \mathcal{F}_t)
\in L^2 (\mathcal{F}_t)$ such that for any $A \in \mathcal{F}_t$ it
holds
$$
\mathcal{E}(\xi 1_A)=\mathcal{E} \left( \mathcal{E}(\xi \vert
\mathcal{F}_t) 1_A \right).
$$
\end{itemize}
\end{definition}

Again in the terminology of \cite{CHMP}, $\mathcal{E}$ is said to
satisfy translation invariance (or to be monetary) if for any $t \in
[0,T]$:
$$
\mathcal{E}(\xi+\eta \vert \mathcal{F}_t)=\mathcal{E}(\xi \vert
\mathcal{F}_t)+\eta, \quad \forall \xi \in L^2(\mathcal{F}_T), \eta
\in L^2(\mathcal{F}_t);
$$
while it is said to be $\mathcal{E}^{\mu}$-dominated (for some $\mu
>0$) if:
$$
\mathcal{E}(\xi+\eta)-\mathcal{E}(\xi) \leq \mathcal{E}^{\mu}(\eta),
\quad \forall \xi, \eta \in L^2(\mathcal{F}_T).
$$

\medskip
\begin{theorem} \label{repr-CHMP}
\textbf{(Coquet et al.; Theorem 7.1; \cite{CHMP})} Let $\mathcal{E}$
be an $\mathcal{F}$-consistent expectation.

If $\mathcal{E}$ satisfies translation invariance and if it is
dominated by some $\mathcal{E} ^{\mu}$ with $\mu >0$, then it is
induced by a conditional $g$-expectation, that is there exists a
function $g$ satisfying (A)-(C), (1${_g}$) such that for any $t \in
[0,T]$
$$
\mathcal{E}(\xi \vert \mathcal{F}_t)=\mathcal{E}_g (\xi \vert
\mathcal{F}_t), \quad \forall \xi \in L^2(\mathcal{F}_t).
$$
\end{theorem}

Some relevant extensions of such a result can be found in Peng
\cite{Peng3} and in Hu et al. \cite{Hu-et-al}, while some applications to risk measures can be found
in Rosazza Gianin \cite{Ros}. The last author, in particular, showed
that $g$-expectations (respectively, conditional $g$-expectations)
provide static (respectively, dynamic) risk measures. More
precisely, the following result holds true. For definitions,
representations and details on (static) risk measures an interested
reader can see Artzner et al. \cite{ADEH}, Delbaen \cite{Delb1},
\cite{Delb2}, F\"{o}llmer and Schied \cite{FS1}, \cite{FS2},
Frittelli and Rosazza Gianin \cite{FrRos1}, among many others.

\begin{proposition} \textbf{(Rosazza Gianin; Proposition 11; \cite{Ros})}
If $g$ satisfies the usual assumptions (including convexity in $z$),
then the risk measure $\rho _g$ defined as
$$
\rho _g (X) \triangleq \mathcal{E}_g (-X)
$$
is a convex risk measure satisfying monotonicity, constancy and
translation invariance.

Moreover: if $g$ also satisfies positive homogeneity in $z$, then
$\rho _g$ is coherent.
\end{proposition}

In view of the result above, some sufficient conditions for a risk
measure to be induced by a $g$-expectation have been found in
\cite{Ros} as an application of Theorem \ref{repr-CHMP}.

\medskip
Note that, at least in the sublinear case and under some suitable
assumptions, one can prove a one-to-one correspondence between the
functional $g$ and the $m$-stable set of generalized scenarios
$\mathcal{S}$ of the suitable risk measure. Hence, one may find (as
an application of the results of Delbaen \cite{Delb3} on $m$-stable
sets) a one-to-one correspondence between time-consistent coherent
risk measures and conditional $g$-expectation. See also Chen and
Epstein \cite{ChEp}.

\medskip
In the sequel, we will prefer to work with concave utilities instead
of convex risk measures. Note that, given a risk measure $\rho$, the
associated monetary utility functional (or, shortly, utility) is
defined as $u \triangleq -\rho$.

\section{Representation of the penalty term of dynamic concave utilities} \label{section3}

In the sequel, we will still work in a Brownian setting, hence
$\mathcal{F}_0$ is trivial. Let $T$ be a fixed finite time horizon.
Given two stopping times $\sigma$ and $\tau $ such that $0 \leq
\sigma \leq \tau \leq T$, consider a \textit{concave monetary
utility} functional $u_{\sigma, \tau}: L^{\infty}
(\mathcal{F}_{\tau}) \to L^{\infty}( \mathcal{F}_{\sigma})$, i.e. a
functional satisfying

\begin{itemize}
\item[(a)] monotonicity: if $\xi,\eta \in L^{\infty}
(\mathcal{F}_{\tau})$ and $\xi \leq \eta$, then
$u_{\sigma,\tau}(\xi) \leq u_{\sigma, \tau}(\eta)$

\item[(b)] translation invariance:
$u_{\sigma,\tau}(\xi+\eta)=u_{\sigma,\tau}(\xi)+\eta$ for any $\xi
\in L^{\infty} (\mathcal{F}_{\tau})$ and $\eta \in L^{\infty}
(\mathcal{F}_{\sigma})$

\item[(c)] concavity: $u_{\sigma,\tau} (\alpha \xi +(1-\alpha)\eta) \geq \alpha u_{\sigma,\tau}(\xi)+ (1-\alpha)
u_{\sigma,\tau}(\eta)$ for any $\xi,\eta \in L^{\infty}
(\mathcal{F}_{\tau})$ and $\alpha \in [0,1]$

\item[(d)] $u_{\sigma,\tau}(0)=0$

\end{itemize}

$(u_{\sigma, \tau}) _{0 \leq \sigma \leq \tau \leq T}$ is called a
\textit{dynamic concave utility}. In particular, $u_{0, T}:
L^{\infty} (\mathcal{F}_T) \to \Bbb R$. The acceptance set
$\mathcal{A}_{\sigma, \tau}$ induced by $u_{\sigma, \tau}$ is
defined as $\mathcal{A}_{\sigma, \tau} \triangleq \{ \xi \in
L^{\infty}( \mathcal{F}_{\tau}) : u_{\sigma, \tau} (\xi) \geq 0 \}$.
To simplify notations, we will often write $u_t$ instead of
$u_{t,T}$.

\bigskip
On $(u_{\sigma, \tau}) _{0 \leq \sigma \leq \tau \leq T}$ we will
assume the following:

\medskip
\noindent \textbf{Assumption (e)}: $(u_{\sigma, \tau}) _{0 \leq
\sigma \leq \tau \leq T}$ is continuous from above (or it satisfies
the Fatou property), i.e. for any decreasing sequence $(\xi_n)_{n
\in \Bbb N}$ in $L^{\infty} (\mathcal{F}_{\tau})$ such that $\lim_n
\xi_n=\xi$ it holds true that $\lim_n u _{\sigma, \tau} (\xi_n)=
u_{\sigma, \tau}(\xi)$.

\medskip
\noindent \textbf{Assumption (f)}: $(u_{\sigma, \tau}) _{\sigma,
\tau}$ is time-consistent, i.e. for all stopping times
$\sigma,\tau,\upsilon$ with $0 \leq \sigma \leq \tau \leq \upsilon
\leq T$:
$$
u_{\sigma,\upsilon}(\xi)=u_{\sigma, \tau} (u_{\tau, \upsilon}(\xi)),
\quad \forall \xi \in L^{\infty}(\mathcal{F}_{\upsilon}).
$$

\medskip
\noindent \textbf{Assumption (g)}: $(u_{\sigma, \tau}) _{\sigma,
\tau}$ satisfies
\begin{equation} \label{000}
u_{\sigma, \tau} (\xi 1_A +\eta 1_{A^c})=u_{\sigma, \tau} (\xi) 1_A
+u_{\sigma, \tau} (\eta) 1_{A^c}, \forall \xi,\eta \in L^{\infty}
(\mathcal{F}_{\tau}), \forall A \in \mathcal{F}_{\sigma} .
\end{equation}

\medskip
\noindent \textbf{Assumption (h)}:  $c_t (P)=0$ for any $t \in
[0,T]$

\medskip It is straightforward to check that this last condition is
equivalent to: $E_P [ \xi \vert \mathcal{F}_t] \geq 0$ for any $\xi
\in \mathcal{A}_t$. Furthermore, $c_0 (P)=0$ can be replaced by the
hypothesis that there is a probability measure $Q$ equivalent to $P$
satisfying $c_0 (Q)=0$.

\medskip
Note that, up to a sign, dynamic concave utilities satisfying the
assumptions above correspond to normalized time-consistent dynamic
risk measures $(\rho_{\sigma, \tau}) _{0 \leq \sigma \leq \tau \leq
T}$ studied, for instance, in Bion-Nadal \cite{BionN} in a general
setting. More precisely, it holds $u_{\sigma, \tau}=-\rho_{\sigma, \tau}$.

\smallskip
By Bion-Nadal \cite{BionN} and Detlefsen and Scandolo \cite{DetlS},
it is known that, under the assumptions above and in the setting of
a general filtration,
\begin{equation} \label{repres}
\begin{array}{ll}
u_{s,t}(\xi)&= \essinf_{Q \sim P, Q =P \mbox{ on } \mathcal{F}_s}
\{E_Q [\xi \vert \mathcal{F}_ {s} ]
+c_{s,t}(Q) \}  \\
&=\essinf_{Q \in \mathcal{P}_{s,t}} \{E_Q [\xi \vert \mathcal{F}_
{s} ] +c_{s,t}(Q) \}
\end{array}
\end{equation}
for any $0 \leq s \leq t \leq T$, where
$$
\begin{array}{rl}
c_{s,t}(Q) &=\esssup _{\xi \in L^{\infty} (\mathcal{F}_t)} \{E_Q
[-\xi \vert \mathcal{F}_s ]+u_{s,t}(\xi) \} \\
\mathcal{P}_{s,t}&=\{Q \mbox{ on } (\Omega, \mathcal{F}_t): Q \sim
P, Q=P \mbox{ on } \mathcal{F}_s \}.
\end{array}
$$
In particular:
$$
\begin{array}{rl}
u_t(\xi) &= \essinf_{Q \sim P; Q =P \mbox{
on } \mathcal{F}_t } \{E_Q [\xi \vert \mathcal{F}_t]+c_{t,T} (Q) \} \\
\\
u_0 (\xi)&=\inf_{Q \sim P}  \{E_Q [\xi ]+c_{0,T} (Q) \}
\end{array}
$$
where $c_t (Q) \triangleq c_{t, T}(Q)=\esssup _{\xi \in
\mathcal{A}_t} E_Q [-\xi \vert \mathcal{F}_t] \geq 0$ and
$\mathcal{A}_t$ denotes the acceptance set induced by $u_t$. Note
that $c_{0,T} (Q)=\sup _{\xi \in L^{\infty}} \{E_Q [-\xi] + u_{0,T}
(\xi) \}$, hence $c_{0,T}$ is lower semi-continuous and is the
Fenchel-Legendre transform of $u$.

Furthermore, Bion-Nadal (see
Theorem 3 in \cite{BionN}) proved that $(\rho_{t,T})_{t \in [0,T]}$
(hence $(u_{t,T})_{t \in [0,T]}$) admits a c\`{a}dl\`{a}g
modification. We will prove in the Appendix that the same is true for $(c_{t,T})_{t \in [0,T]}$.

Note that in \cite{BionN} and \cite{DetlS} the representation
(\ref{repres}) was shown with $Q \ll P$ instead of $Q \sim P$.
Nevertheless, assumption (h) guarantees that the representation (\ref{repres}) also holds true (for a proof see Kl\"{o}ppel
and Schweizer \cite{KlS} and, in discrete-time, Cheridito et al.
\cite{CDK3} and F\"{o}llmer and Penner \cite{FP}).

\begin{remark} \label{rem1}
It is evident that if $(u_t)_{t\geq 0}$ is time-consistent, if
$u_t(0)=0$ and if it satisfies condition (\ref{000}), then
$$
u_0 (\xi 1_A)=u_0(u_t(\xi1_A))=u_0 (u_t(\xi) 1_A)
$$
for any $t \geq 0$, $\xi \in L^{\infty} (\mathcal{F}_T)$ and $A \in
\mathcal{F}_t$.
\end{remark}

It is therefore clear that if $(u_{\sigma, \tau}) _{\sigma, \tau}$
is time-consistent, then everything is defined by $u_0$. The
relevance of time-consistency of the dynamic concave utility is also
underlined by the following results. On one hand, as shown by
Delbaen \cite{Delb3} and Cheridito et al. \cite{CDK3},
time-consistency is indeed equivalent to the \textit{decomposition
property} of acceptable sets, that is
$$
\mathcal{A}_{\sigma,
\upsilon}=\mathcal{A}_{\sigma,\tau}+\mathcal{A}_{\tau, \upsilon}
$$
for all stopping times $\sigma, \tau, \upsilon$ such that $0 \leq
\sigma\leq\tau\leq\upsilon\leq T$. On the other hand, both the
properties above are equivalent to the \textit{cocycle property} of
the penalty term $c$, that is
$$c_{\sigma,\upsilon} (Q)=c_{\sigma,\tau} (Q)+E_Q [c_{\tau,\upsilon}
(Q) \vert \mathcal{F}_{\sigma}]$$ for all stopping times $\sigma,
\tau, \upsilon$ such that $0 \leq \sigma \leq \tau \leq \upsilon
\leq T$ (see Bion-Nadal \cite{BionN} for the definition and the
proof).

\bigskip
In the sequel, we use the terminology of Rockafellar \cite{Rock2} on
convex functions. Our aim is now to prove the following result.

\begin{theorem} \label{thm3} Let $(u_{\sigma, \tau}) _{0 \leq \sigma \leq \tau \leq
T}$ be a dynamic concave utility satisfying the assumptions above.

\begin{itemize}
\item[(i)] For all stopping times $\sigma, \tau$ such that $0 \leq \sigma \leq \tau \leq
T$ and for any probability measure $Q$ equivalent to $P$:
\begin{equation} \label{c-f}
c_{\sigma, \tau} (Q)=E_Q \left[ \int _{\sigma} ^{\tau} f(u,q_u)
du | \mathcal{F}_{\sigma} \right]
\end{equation}
for some suitable function $f:[0,T] \times \Omega \times \Bbb R ^d
\to [0, + \infty]$ such that $f(t, \omega, \cdot)$ is proper, convex
and lower semi-continuous.

\item[(ii)] For all stopping times $\sigma, \tau$ such that $0 \leq \sigma \leq \tau \leq T$ and $\xi \in
L^{\infty} (\mathcal{F}_T)$, the dynamic concave utility in
(\ref{repres}) can be represented as
$$
u_{\sigma,\tau}(\xi)=\essinf _{Q \in \mathcal{P}_{\sigma,\tau}} E_Q
\left[ \xi+ \int _{\sigma} ^{\tau} f(u,q_u) du |
\mathcal{F}_{\sigma} \right].
$$
\end{itemize}
\end{theorem}

\medskip
\begin{remark}
For a dynamic concave utilities satisfying assumptions (e), (g), (h), from Theorem 1 of Bion-Nadal \cite{BionN} it follows that
Theorem \ref{thm3}(i) is equivalent to time-consistency (assumption
(f)) of  $(u_{\sigma, \tau}) _{0 \leq
\sigma \leq \tau \leq T}$.
\end{remark}

\begin{remark}
In an incomplete market, the lower price $\inf _{Q \in \mathcal{M}}
E_Q [\xi]$ (where $\mathcal{M}$ denotes the set of all risk-neutral
probability measures) defines a utility satisfying all our
properties but it is not given by a $g$-expectation. See Delbaen
\cite{Delb3} for details about how to get $f$.
\end{remark}

The proof of Theorem \ref{thm3} will be decomposed into several
steps as outlined below.

Set
\begin{equation} \label{u^n}
u_{s,t} ^n(\xi)=\essinf _{Q \sim P; \Vert q\Vert \leq n} \{E_Q [\xi
\vert \mathcal{F}_s ]+c_{s,t} (Q) \}.
\end{equation}

\smallskip
Note that (by definition of $u^n$ and by assumption (h)) for any
$\xi \in L^{\infty}(\mathcal{F}_T)$ it holds $u_t^0 (\xi)=E_P [\xi
\vert \mathcal{F}_t]$ and $u^n _t (\xi)\leq E_P [\xi \vert
\mathcal{F}_t]$.

\medskip
\begin{remark} \label{remark-wcomp} The reason why the truncated
utility $u^n$ has been defined as above is due to the fact that the
set $\{Q \sim P; \Vert q\Vert \leq n \}$ is weakly compact. This
argument will be useful in the proof of Proposition \ref{prop-f}.
\end{remark}

\bigskip
\begin{proposition} \label{prop-f}
Suppose that the dynamic concave utility $(u_{\sigma, \tau}) _{0
\leq \sigma \leq \tau \leq T}$ satisfies the assumptions above.
Then:

\begin{enumerate}
\item[(i)] $u^n$ is a dynamic concave utility satisfying assumptions (e)-(g). Moreover,
the acceptance sets induced by $u^n$ satisfy the decomposition property and
\begin{equation} \label{c_n}
c^n_{s,t} (Q)= \left\{
\begin{array} {ll}
c_{s,t}(Q); \quad &\mbox{ if } \quad \Vert q \Vert \leq n \\
+ \infty ; \quad &\mbox { otherwise }
\end{array}
\right.
\end{equation}
satisfies the cocycle property and $c^n_{s,t} (P)=0$.

\item[(ii)] $u^n$ is induced by a conditional $g_n$-expectation, i.e.
$$
u^n_t(\xi)=- \mathcal{E}_{g_n} (-\xi \vert \mathcal{F}_t)
$$
for some convex function $g_n$ satisfying the usual conditions and
such that $g_n(\cdot, \cdot, z)$ is predictable for any $z \in \Bbb
R ^d$. In other words, $u^n$ satisfies the following BSDE
\begin{equation} \label{u}
\left\{
\begin{array}{rl}
du^n_t(\xi)&=g_n(t,Z^n_t)dt-Z^n_t dB_t \\
u^n_T(\xi)&=\xi
\end{array}
\right.
\end{equation}

\item[(iii)] For any probability measure $Q \sim P$ such that $\Vert q \Vert \leq
n$ it holds that for any $0 \leq s \leq t \leq T$:
$$
\begin{array}{rl}
c_{0,t} ^n (Q)&= E_Q \left[ \int _0 ^{t} f_n (u,q_u) du \right] \\
c_{s,t} ^n (Q)&= E_Q \left[ \int _s ^{t} f_n (u,q_u) du |
\mathcal{F}_s \right]
\end{array}
$$
where $f_n: [0,T] \times \Omega \times \Bbb R^d \to [0, +\infty]$ is
induced (by duality) by $g_n$ and $f_n (t, \omega, \cdot)$ is
proper, convex and lower semi-continuous.

\item[(iv)] The sequence of convex functions $g_n$ is increasing in $n$.

\item[(v)] The sequence of $f_n$ is decreasing in $n$ and, for any $n \geq 0$,
$f_n (t, \omega,q)=+\infty$ for $\Vert q \Vert >n$.

Furthermore, once $(t, \omega)$ is fixed, for any $q$ either there
exists $n \geq 0$ such that
$$f_n (t,\omega, q)=f_m(t, \omega, q)= f(t,\omega,
q) < +\infty, \quad \forall m \geq n$$ or for all $n \geq 0$
$$f_n (t,\omega, q)=+\infty=f(t, \omega, q),$$
for some function $f: [0,T] \times \Omega \times \Bbb R^d \to [0,
+\infty]$.

Hence $f(t, \omega, x)=\inf _n f_n (t, \omega, x)$ and it is such
that $f(t, \omega, \cdot)$ is proper, convex and lower
semi-continuous.
\end{enumerate}

\end{proposition}

\begin{proof}
\begin{enumerate}

\item[(i)]
From the representation (\ref{u^n}) it follows that $u^n$ is a
dynamic concave utility which is continuous from above (see
Detlefsen and Scandolo \cite{DetlS} and Kl\"{o}ppel and Schweizer
\cite{KlS}). Still from (\ref{u^n}) one deduces that $u_{\sigma,
\tau}^n (\xi1_A)=u_{\sigma, \tau} ^n (\xi) 1_A$ for any $\xi \in
L^{\infty}(\mathcal{F}_T)$, $0 \leq \sigma \leq \tau \leq T$ and $A
\in \mathcal{F}_{\sigma}$. Hence, by Proposition 2.9 of Detlefsen
and Scandolo \cite{DetlS}, also assumption (g) is satisfied.

The cocycle property of $c^n$ and time-consistency of $u^n$ follow
from
$$
c_{s,t}^n (Q)= \left \{ \begin{array}{ll} c_{s,t}(Q); \quad &\mbox{
if } \quad \Vert q\Vert \leq n \\
+ \infty; \quad &\mbox{ otherwise}
\end{array} \right.
$$
and from Theorem 1 of Bion-Nadal \cite{BionN}.

Since for the probability measure $P$ it holds $q^P \equiv 0$,
$c^n_{s,t} (P)=c_{s,t} (P)=0$.

The decomposition property of acceptance sets is due to Theorem
4.6 of Cheridito et al. \cite{CDK3} and, later, to Theorem 1 of
Bion-Nadal \cite{BionN}.

\item[(ii)] Set $\pi _t ^n (\xi) \triangleq -u_t ^n (-\xi)=\esssup _{Q \sim P; \Vert q\Vert \leq n} \{E_Q [\xi
\vert \mathcal{F}_t ]-c_{t} (Q) \}$.

From (i), $(\pi ^n _{\sigma, \tau})_{0 \leq \sigma \leq \tau \leq
T}$ is time-consistent. Furthermore, it is easy to check that it
satisfies monotonicity, translation invariance and constancy (this
last follows from the assumption $c_t(P)=0$).

Moreover, $\pi _0 ^n$ satisfies strict monotonicity. This property
follows from weak compactness of the set $\{Q \sim P: \Vert q \Vert
\leq n \}$ (see Remark \ref{remark-wcomp}). In order to verify
strict monotonicity, consider $\eta \geq \xi$ such that
$P(\eta>\xi)>0$. Since $\pi _0 ^n (\xi)=E_Q[\xi]-c_0 (Q)$ for some
$Q \sim P$ such that $\Vert q \Vert \leq n$, $\pi _0 ^n (\eta) \geq
E_Q [\eta]-c_0 (Q) > E_Q [\xi]-c_0 (Q)=\pi _0 ^n (\xi)$.

Finally, we will show that $\pi _0 ^n$ is dominated by some
$\mathcal{E}^{\mu}$. For any $\xi,\eta \in
L^{\infty}(\mathcal{F}_T)$
$$
\begin{array}{rl}
&\pi _0 ^n (\xi+\eta)-\pi _0 ^n (\xi) \\
&=\sup _{Q: \Vert q \Vert \leq n} \{ E_Q[\xi+\eta]-c_0 (Q) \}-\sup
_{Q: \Vert q \Vert \leq n} \{ E_Q[\xi]-c_0
(Q) \} \\
&\leq \sup _{Q: \Vert q \Vert \leq n}  E_Q [\eta]= \mathcal{E}^n
(\eta).
\end{array}
$$
The last equality follows from Lemma 3 of Chen and Peng \cite{ChP}
($\Bbb R$ case) which may be extended to $\Bbb R ^d$.

By the arguments above and Remark \ref{rem1}, $(\pi _t ^n)_{t \geq
0}$ satisfies the hypothesis of Theorem \ref{repr-CHMP}. Hence there
exists a functional $g_n:[0,T] \times \Omega \times \Bbb R ^d \to
\Bbb R$ satisfying assumptions (A)-(C), (1$_g$) and such that $ \pi
_t ^n (\xi)= \mathcal{E} _{g_n} (\xi \vert \mathcal{F}_t) $.

\smallskip
It can be checked that $g_n(\cdot, \cdot, z)$ is predictable for any
$z \in \Bbb R^d$ (see also Theorem 3.1 of Peng \cite{Peng3}).
Furthermore, since $\pi _t ^n$ is a convex functional, by Theorem
3.2 of Jiang \cite{Jiang} it follows that $g_n(t, \omega, \cdot)$
has to be convex. Hence
$$
\begin{array}{rl}
u_t ^n (\xi)&=-\mathcal{E} _{g_n} (-\xi \vert \mathcal{F}_t) \\
u_0 ^n (\xi)&=-\mathcal{E} _{g_n} (-\xi)
\end{array}
$$
for some function $g_n$ satisfying the usual conditions. It is
therefore immediate to check that $u^n$ satisfies the BSDE in
(\ref{u}).

Moreover: for almost all $(t, \omega)$ it holds that the set $\{z
\in \Bbb R ^d : g_n(t, \omega, z) \leq \alpha \}$ is closed for any
$\alpha \in \Bbb R$. The closure of such a set (or, equivalently,
the lower semi-continuity of $g_n (t, \omega, \cdot)$) is due to the
fact that $g_n$ is Lipschitz with constant $n$ (see the arguments
above and Theorem \ref{repr-CHMP}). Hence $g_n(t, \omega, \cdot)$ is
convex, proper and lower semi-continuous.

\item[(iii)] Set now
\begin{equation} \label{fn} f_n(t, \omega, q) \triangleq
\sup_{z \in \Bbb R ^d} \{q \cdot z -g_n(t, \omega, z)\}.
\end{equation}
Note that $f_n(t, \omega,q) \geq 0$ (take for instance $z =0$ in the
definition of $f_n$) and, because of the assumption $c_0 (P)=0$,
$f_n(t,0)=0$.

Since $g_n (t, \omega,z)$ is predictable (by item (ii)),
$$
f_n(t, \omega, q) = \sup_{z \in \Bbb R ^d} \{q \cdot z -g_n(t,
\omega, z)\}= \sup_{z \in \Bbb Q ^d} \{q \cdot z -g_n(t, \omega,
z)\}
$$
is predictable for any $q \in \Bbb R ^d$ (as supremum of countably
many predictable elements). Note that $\Vert q \Vert > n$ implies
$f_n (t, \omega, q)=+ \infty$ (by (\ref{fn})).

\medskip
Since $g_n(t, \omega, \cdot)$ is convex, proper and lower
semi-continuous (see above) and $f_n(t, \omega, \cdot)$ is the
convex conjugate of $g_n(t, \omega, \cdot)$, i.e. $f_n( q)=g_n ^*
(q)$, also $f_n$ is convex, proper and lower semi-continuous (see
Rockafellar \cite{Rock2}).

\smallskip
As a consequence of the dual representation of a $g$-expectation in
Theorem 7.4 of Barrieu and El Karoui \cite{BaElK2} we get
$$c_{0,T} ^n (Q)=E_Q \left[\int_0 ^{T} f_n(u,q_u) du \right]$$ for any probability measure $Q \sim P$ such
that $\Vert q \Vert \leq n$.

\medskip
It remains to show that $c_{s,t} ^n (Q)=E_Q \left[\int_s ^{t}
f_n(u,q_u) du | \mathcal{F}_s \right]$ for any $0 \leq s \leq
t \leq T$ and for any probability measure $Q \sim P$ such that
$\Vert q \Vert \leq n$. Also this result can be deduced by Theorem
7.4 of Barrieu and El Karoui \cite{BaElK2}. Nevertheless, since the
proof will be useful later, we postpone it to Lemma
\ref{lemma-proof}.

\item[(iv)] It is easy to check that the sequences of $u_0 ^n$ and of $c_0 ^n$ are decreasing in $n
\in \Bbb N$. By applying the Converse Comparison Theorem on BSDE
(see Briand et al. \cite{BCHMP}) and Lemma 2.1 of Jiang
\cite{Jiang}, we will show that the sequence of convex functions
$g_n$ (which induce $u^n$) is increasing in $n$.

In order to prove the thesis above we will proceed in a similar way
as in Jiang \cite{Jiang}. By definition of $u^n$, $u_{0,T} ^n(\xi)
\geq u_{0,T} ^{n+1} (\xi)$ as well as $u_{s,T} ^n (\xi) \geq u_{s,T}
^{n+1} (\xi)$ hold true for any $\xi \in L^{\infty}
(\mathcal{F}_T)$. By item (ii) we deduce therefore that for any $\xi
\in L^{\infty} (\mathcal{F}_T)$
\begin{equation} \label{ineqE}
\begin{array}{rl}
\mathcal{E} _{g_n} (\xi) &\leq \mathcal{E} _{g_{n+1}} (\xi) \\
\mathcal{E} _{g_n} (\xi \vert \mathcal{F}_s )&\leq \mathcal{E}
_{g_{n+1}} (\xi \vert \mathcal{F}_s)
\end{array}
\end{equation}
Denote now by $\mathcal{E}_g ^{s,t}$ the conditional $g$-expectation
at time $s$ with final time $t$. To apply successfully Lemma 2.1 of
Jiang \cite{Jiang} we need to verify that
\begin{equation} \label{conditionE}
\mathcal{E}_{g_n} ^{s,t}(\xi) \leq \mathcal{E}_{g_{n+1}} ^{s,t}
(\xi), \quad \forall s,t \in [0,T] \mbox{ with } s \leq t, \quad
\forall \xi \in L^{\infty} (\mathcal{F}_t).
\end{equation}

Condition (\ref{conditionE}) has already been established for
$(s,t)=(0,T)$ and $(s,t)=(s,T)$. Consider now the case
$(s,t)=(0,t)$. Since (see Peng \cite{Peng1} for details)
$\mathcal{E}_{g} ^{s,t}(\eta) = \mathcal{E}_{g} ^{s,T} (\eta)$ for
any $\eta \in L^{\infty} (\mathcal{F}_t)$, from (\ref{ineqE}) we
deduce that $\mathcal{E}_{g_n} ^{0,t}(\xi) \leq
\mathcal{E}_{g_{n+1}} ^{0,t} (\xi)$ for any $0\leq t \leq T$ and
$\xi \in L^{\infty} (\mathcal{F}_t)$. For general $(s,t)$,
inequality (\ref{ineqE}) can be checked as above. Indeed, for any
$\xi \in L^{\infty} (\mathcal{F}_t)$ it holds that
$\mathcal{E}_{g_n} ^{s,t}(\xi)=\mathcal{E}_{g_n} ^{s,T}(\xi) \leq
\mathcal{E}_{g_{n+1}} ^{s,T} (\xi)=\mathcal{E}_{g_{n+1}}
^{s,t}(\xi)$.

Set now
$$
\mathcal{S} ^z (g) \triangleq \{ t \in [0,T): g(t,z)= L^1-\lim
_{\varepsilon \to 0 ^+ } \frac{1}{\varepsilon} \mathcal{E}_g ^{t, t+
\varepsilon} (z(B_{t+\varepsilon} -B_t) )\}.
$$

From Lemma 2.1 of Jiang \cite{Jiang} it follows that
$$
m([0,T) \setminus \mathcal{S}^z (g_i))=0 \quad \forall z \in \Bbb R
^d
$$
for $i=n,n+1$, where $m$ denotes the Lebesgue measure on $[0,T]$.

By the arguments above it follows that for any $z \in \Bbb R ^d$
$$
\mbox{ if } t \in \mathcal{S}^z (g_n) \cap \mathcal{S}^z (g_{n+1})
\neq \emptyset \qquad \Rightarrow g_n(t,z) \leq g_{n+1} (t,z) \quad
P\mbox{-a.s.}
$$
and
$$
\begin{array}{rl}
&m([0,T) \setminus (\mathcal{S}^z (g_n)\cap \mathcal{S}^z
(g_{n+1}))) \\
&= m(([0,T) \setminus \mathcal{S}^z (g_n)) \cup ([0,T) \setminus
\mathcal{S}^z (g_{n+1})))=0
\end{array}
$$
Hence, by proceeding as in Jiang \cite{Jiang} it can be checked that
for any $z \in \Bbb R ^d$
$$
g_n(t,z) \leq g_{n+1} (t,z) \quad (dt \times dP)\mbox{-a.s.}
$$

Positivity of any $g_n$ is due to the fact that $u_t ^0 (\xi)=E_P
[\xi \vert \mathcal{F}_t ]=- \mathcal{E} _{g_0} (-\xi \vert
\mathcal{F}_t)$ where $g_0 \equiv 0$. By the same arguments above,
therefore, $g_n \geq g_0 \equiv 0$.

\item[(v)]  From (iii) and (iv) it then follows that the sequence of $f_n$ is decreasing in
$n$.

Consider again the measurable space $([0,T] \times \Omega,
\mathcal{P}, m \times P)$, where $\mathcal{P}$ denotes the
predictable $\sigma$-algebra and $m$ denotes the Lebesgue measure on
$[0,T]$. Denote by $\overline{\mathcal{P}}$ the completion of
$\mathcal{P}$.

Take $N>0$ and, for any $\varepsilon>0$, set
$$
\begin{array}{rl}
&E=E_{N,\varepsilon} \\
&\triangleq \left\{(t, \omega, q) \in [0,T] \times \Omega \times
\Bbb R ^d \vert \begin{array}{c}
\mbox{ } \Vert q \Vert \leq n; \\
f_{n+1} (t, \omega, q) + \varepsilon <f_n ( t, \omega, q) \leq N
\end{array} \right\}
\end{array}
$$
and $\pi(E)$ its projection on $[0,T] \times \Omega$. Note that $E
\in \overline{\mathcal{P}} \otimes \mathcal{B} (\Bbb R ^d)$.

From the Measurable Selection Theorem (see Aumann \cite{Aumann} and
Aliprantis and Border \cite{AB}), $\pi (E) \in
\overline{\mathcal{P}}$ and there exists a $\overline{\mathcal{P}}$-
measurable $\overline{q}: \pi (E) \to \Bbb R ^d$ such that $(t,
\omega, \overline{q}(t, \omega) ) \in E$ for $(m \times P)$-a.e.
$(t, \omega) \in \pi(E)$. Set now $\overline{q}=0$ on $\pi(E) ^c$.
To such a $\overline{q}$ it is therefore possible to associate a $q:
[0,T] \times \Omega \to \Bbb R ^d$ which is $\mathcal{P}$-measurable
and equal to $\overline{q}$ $(m \times P)$-almost everywhere.

Let $Q$ be the probability measure associated to $q$ as above. By
definition, $\Vert q \Vert \leq n$. Hence, $c_{0,T} ^n (Q)=c_{0,T}
^{n+1} (Q)=c_{0,T} (Q) <+ \infty $. Furthermore, by definition of
$E$ it follows that
$$
\begin{array}{rl}
c_{0,T} ^n (Q)&= E_Q \left[\int_0 ^T f_n (u,q_u) du \right] \\
&=E_Q \left[\int_0 ^T f_n (u,q_u) 1_{\pi (E)} du \right]+ E_Q
\left[\int_0 ^T f_n (u,q_u) 1_{\pi (E) ^c} du \right] \\
&=E_Q \left[\int_0 ^T f_n (u,q_u) 1_{\pi (E)} du \right]  \\
& \geq E_Q \left[\int_0 ^T [f_{n+1} (u,q_u)+ \varepsilon] 1_{\pi
(E)}
du \right] \\
&=E_Q \left[\int_0 ^T f_{n+1} (u,q_u) 1_{\pi (E)} du \right]
+\varepsilon (m \times Q)( \pi (E) ) \\
&=c_{0,T} ^{n+1} (Q) + \varepsilon (m \times Q ) ( \pi (E))
\end{array}
$$
If  $(m \times Q)(\pi (E))>0$, then $c_{0,T} ^{n+1} (Q) +
\tilde{\varepsilon} < c_{0,T} ^n (Q) < +\infty$, that is a
contradiction. Hence, $(m \times Q)(\pi (E))=0$, i.e. $(m \times
Q)(\{(t, \omega) : N \geq f_n (t,\omega, q_t) > f_{n+1} (t,\omega,
q_t )+ \varepsilon \})=0$.

By letting $N$ tend to $+\infty$, from the arguments above and since
$Q \sim P$, it follows that if $f_n < + \infty$ on $\{x: \Vert x
\Vert \leq n \}$
$$
f_n = f_{n+1} \quad (m \times dP)\mbox{-a.s.}
$$
and hence $f_n = f_{n+1}=f$ $(m \times dP)$-a.s. for some functional
$f$. I.e.
$$
f_n(t, \omega, x)=f_{n+1} (t, \omega, x)=f(t, \omega, x)  \quad (m
\times dP) \mbox{-a.s.}, \quad \mbox{ for } \Vert x \Vert \leq n .
$$

Furthermore, we may conclude that, once $(t, \omega)$ is fixed,  for
any $q$ either (1) there exists $n \geq 0$ such that $f_n(t, \omega,
q)< + \infty $ (hence $f_m(t, \omega, q)=f(t, \omega, q)< +\infty$
for any $m \geq n$ and $m \geq \Vert q \Vert$) or (2) for all $n\geq
0$ it holds $f_n(t, \omega, q)=+\infty=f(t, \omega, q)$. Hence
$$
f(t, \omega, x)=\inf_n f_n (t, \omega, x).
$$

\medskip
By the properties of the sequence of $f_n$, it follows that
$f(\cdot, \cdot,0)=0$.

It remains to prove that $f(t, \omega, \cdot)$ is proper, convex and
lower semi-continuous. Properness of $f(t, \omega, \cdot)$ is
trivial. Since $f(t, \omega, x)=\lim _n f_n (t, \omega, x)=\inf _n
f_n (t, \omega, x)$ for almost all $(t, \omega)$ and any $f_n$ is
predictable and convex in $x$, it is easy to check that also $f$ is
predictable and convex in $x$.

Furthermore, for almost all $(t, \omega)$ the set $\{q \in \Bbb R ^d
: f(t,\omega, q) \leq \alpha \}$ is closed for any $\alpha \in \Bbb
R$. Take indeed a sequence $\{q^k \}_{k \geq 0}$ such that $q^k \to
_k q$ and $f(t, \omega, q^k) \leq \alpha$. There exists $N \in \Bbb
N$ such that $\Vert q^k \Vert \leq N$ for all $k$. Hence
$$
f(t, \omega, q^k)=f_N(t, \omega, q^k)\leq \alpha .
$$
Since $f_N(t, \omega, \cdot)$ is lower semi-continuous, $f(t,
\omega, q)=f_N(t, \omega, q) \leq \underline{\lim}_k f_N (t, \omega,
q^k) \leq \alpha$. Hence also $f(t, \omega, \cdot)$ is lower
semi-continuous.
\end{enumerate}
\end{proof}

\begin{lemma}\label{lemma-proof}
If $c_{0,T} ^n (Q)=E_Q \left[\int_0 ^{T} f_n(u,q_u) du \right]$
holds for any probability measure $Q \sim P$ such that $\Vert q
\Vert \leq n$, then also $c_{s,t} ^n (Q)=E_Q \left[\int_s ^{t}
f_n(u,q_u) du | \mathcal{F}_s \right]$ holds for any $0 \leq s
\leq t \leq T$ and for any probability measure $Q \sim P$ such that
$\Vert q \Vert \leq n$.
\end{lemma}

\begin{proof}
Let $Q$ be a probability measure equivalent to $P$ and such that
$\Vert q \Vert \leq n$. Consider the case where $s=0$ and take the
probability measure $\overline{Q}$ corresponding to the following
$\overline{q}$:
$$
\overline{q}_u= \left\{ \begin{array}{lll} &q_u; &\mbox{ if } 0 \leq
u
\leq t \\
&0; &\mbox{ if } t < u \leq T \end{array} \right.
$$
obtained by pasting $Q$ and $P$. It is clear that $\Vert
\overline{q} \Vert \leq n$. From the cocycle property of $c^n$
established in Proposition \ref{prop-f}(i) it follows that
$$
\begin{array}{rl}
c_{0,T} ^n (\overline{Q})&=c_{0,t} ^n
(\overline{Q})+E_{\overline{Q}} [c_{t,T}^n (\overline{Q})] \\
&=c_{0,t} ^n (\overline{Q})+E_{\overline{Q}} [c_{t,T}^n (P)]=c_{0,t}
^n (\overline{Q}).
\end{array}
$$
From the arguments above it follows that
$$
\begin{array}{rl}
c_{0,t} ^n (Q)&=c_{0,t} ^n (\overline{Q})=c_{0,T} ^n (\overline{Q})
\\
&=E_{\overline{Q}} \left[ \int _0 ^T f_n(u, \overline{q}_u) du
\right] \\
&=E_{\overline{Q}} \left[ \int _0 ^t f_n(u, \overline{q}_u) du
\right]=E_{Q} \left[ \int _0 ^t f_n(u, q_u) du \right].
\end{array}
$$

We will now come back to the general case. Consider the probability
measure $Q^*$ obtained by pasting $Q$ and $P$ as follows:
$$
q^*_u= \left\{ \begin{array}{lll} &0; &\mbox{ if } 0 \leq u
\leq s \\
&q_u 1_A +0 \cdot 1_{A^c} ; &\mbox{ if } s<u \leq T \end{array}
\right. ,
$$
with $A \in \mathcal{F}_s$. On one hand, we deduce that
$c_{0,s}^n(Q^*)=c_{0,s}^n (P)=0$, while for any $s <t \leq T$
$$
\begin{array}{rl}
c_{0,t} ^n (Q^*)&= E_{Q^*} \left[\int _0 ^t f_n (u, q_u^*)du \right]
\\
&=E_{Q^*} \left[1_A\int _s ^t f_n (u, q _u) du  \right] \\
&=E_P \left[ E_Q \left[1_A \int _s ^t f_n (u,q_u) du  |
\mathcal{F}_s \right] \right] \\
&=E_P \left[1_A E_Q \left[\int _s ^t f_n (u,q_u) du  |
\mathcal{F}_s \right] \right].
\end{array}
$$
On the other hand, from the cocycle property $E_{Q^*} [c_{s,t} ^n
(Q^*)]=c_{0,t}^n (Q^*) -c_{0,s} ^n (Q^*)$, hence
$$
\begin{array}{rl}
c_{0,t}^n (Q^*)&=c_{0,t}^n (Q^*) -c_{0,s} ^n (Q^*)=E_{Q^*} [c_{s,t}
^n (Q^*)] \\
&=E_{Q^*}[E_{Q^*} [c_{s,t} ^n (Q^*) \vert \mathcal{F}_s]] \\
&=E_{P}[1_A E_{Q} [c_{s,t} ^n (Q) \vert \mathcal{F}_s]].
\end{array}
$$
Since the set $A$ is arbitrary, we deduce that for any $A \in
\mathcal{F}_s$
$$
E_P \left[1_A E_Q \left[\int _s ^t f_n (u,q_u) du  |
\mathcal{F}_s \right] \right]=E_{P} [1_A E_{Q} [c_{s,t} ^n (Q) \vert
\mathcal{F}_s ]],
$$
hence
$$
c_{s,t}^n(Q)=E_Q[c_{s,t}^n(Q) \vert \mathcal{F}_s]=E_Q
\left[\int _s ^t f_n (u,q_u) du | \mathcal{F}_s \right].
$$
\end{proof}

\begin{lemma} \label{lemma4}
For any probability measure $Q$ equivalent to $P$ it holds true that
$$
\begin{array}{ll}
c_{0,T} (Q) &\leq E_Q \left[ \int _0 ^T f(u,q_u) du \right] \\
c_{t,T} (Q) &\leq E_Q \left[ \int _t ^T f(u,q_u) du |
\mathcal{F} _t \right]
\end{array}
$$
\end{lemma}

\begin{proof} We will start proving the inequality for $c_{0,T}(Q)$.

\medskip
\noindent \textit{Case 1:} $\int _0 ^T f(u,q_u) du $ is bounded.

Consider the probability measure $Q^n$ corresponding to $q^{n}
\triangleq q 1_{\Vert q \Vert \leq n}$. Since $\int _0 ^T f(u,q_u)
du $ is bounded (by assumption), from Proposition \ref{prop-f}(iii)
it follows that
$$
\begin{array}{rl}
\lim _n c_{0,T}(Q^n)&= \lim _n E_{Q^n} \left[\int _0 ^T f_n (u, q_u
^n)
du \right] \\
&=\lim _n E_{Q^n} \left[\int _0 ^T f(u, q_u ) 1_{\Vert q \Vert \leq
n} du
\right]\\
&=\lim _n E_{Q} \left[\frac{dQ^n}{dQ} \int _0 ^T f(u, q_u ) 1_{\Vert
q \Vert \leq n} du \right] \\
&=E_Q \left[\int _0 ^T f(u, q_u ) du \right] < + \infty
\end{array}
$$
Since $\frac{dQ ^{n}}{dP} \to ^{L^1} _n \frac{dQ}{dP}$, by lower
semi-continuity of $c_{0,T}(Q)$ it follows that
$$c_{0,T}(Q) \leq \liminf _n c_{0,T}(Q^n) \leq E_Q \left[\int _0 ^T f(u, q_u ) du
\right].$$

\medskip
\noindent \textit{Case 2:} $\int _0 ^T f(u,q_u) du \in L^1 (Q)$.

For any $n \in \Bbb{N}$, set $\sigma _n \triangleq \inf \{t \geq 0 :
\int _0 ^t f(u,q_u)du \geq n \}$. Then $\sigma_n$ is a stopping time
and $\sigma _n \uparrow T$.

Set $Q^{\sigma _n}$ the probability measure corresponding to
$\frac{dQ^{\sigma _n}}{dP}=\mathcal{E} (q \cdot B) ^{\sigma _n}$. It
is easy to check that $\frac{dQ ^{\sigma _n}}{dP} \to ^{L^1} _n
\frac{dQ}{dP}$. Furthermore,
$$
E_{Q^{\sigma_n}} \left[ \int _0 ^{\sigma _n} f(u,q_u^{\sigma _n}) du
\right] =E_{Q} \left[ \int _0 ^{\sigma _n} f(u,q_u) du \right] \to
_n E_Q \left[ \int _0 ^{T} f(u,q_u) du \right],
$$
where the equality above is due to the fact that $q$ and
$q^{\sigma_n}$ coincide on the stochastic interval $\sintl 0, \sigma
_n \sintr$. By applying the arguments above, we obtain
$$
\begin{array}{ll}
c_{0,T}(Q) &\leq \liminf _n c_{0,T}(Q ^{\sigma _n}) \leq \liminf _n
E_{Q^{\sigma_n}} \left[ \int _0 ^{\sigma _n} f(u,q_u^{\sigma _n}) du
\right] \\
&\leq E_Q \left[ \int _0 ^{T} f(u,q_u) du \right].
\end{array}
$$

\medskip
\noindent \textit{Case 3:} General case.

In general, if $\int _0 ^T f(u,q_u) du \notin L^1 (Q)$, then $E_Q
\left[ \int _0 ^{T} f(u,q_u) du \right]=+\infty$. Hence
$c_{0,T}(Q)\leq E_Q \left[ \int _0 ^{T} f(u,q_u) du \right]$.

\medskip
The inequality $c_{t,T} (Q) \leq E_Q \left[ \int _t ^T f(u,q_u)
du | \mathcal{F}_t \right]$ can be checked by proceeding as in
the proof of Lemma \ref{lemma-proof}.
\end{proof}

\begin{lemma} \label{lemma6}
Let $Q$ be a probability measure equivalent to $P$ and such that
$c_{0,T}(Q)<+\infty$. If $\{ \tau _n \}_{n \geq 0}$ is a sequence
of stopping times such that $P(\tau _n <T) \to _n 0$, then $c_0
(Q^{\tau _n}) \uparrow c_0 (Q)$, where $Q^{\tau _n}$ is defined
by $\frac{dQ ^{\tau _n}}{dP}=E_P \left[\frac{dQ}{dP} \vert
\mathcal{F}_{\tau _n} \right]$.
\end{lemma}

\begin{proof}
On one hand, by the cocycle property and by the definition of
$Q^{\tau _n}$ it follows
$$
\begin{array}{rl}
c_{0,T} (Q)&=c_{0, \tau _n} (Q)+E_Q [c_{\tau _n, T} (Q)] \\
&\geq c_{0, \tau _n} (Q) =c_{0,T} (Q^{\tau _n})
\end{array}
$$
On the other hand, by the lower semi-continuity of $c_0$ and by $\frac{dQ ^{\tau _n}}{dP} \to ^{L^1} _n \frac{dQ}{dP}$ it holds
$c_{0,T} (Q) \leq \liminf _n c_{0,T} (Q^{\tau _n})$. So $\lim _n
c_{0,T} (Q^{\tau _n})=c_{0,T}(Q)$.
\end{proof}

\begin{lemma} \label{lemma7} Consider a general setting where the
filtration satisfies the usual hypothesis but it is not necessarily
a Brownian filtration.

Let $Q$ be a probability measure equivalent to $P$ such that
$c_{0,T}(Q)<+\infty$ and $(c_t(Q))_{t \in [0,T]}$ is
right-continuous.

Then there exists a unique increasing, predictable process $(A_t)_{t
\in [0,T]}$ (depending on $Q$) such that $A_0=0$ and
\begin{equation} \label{c-A}
c_t(Q)=E_Q [A_T - A_t \vert \mathcal{F}_t ], \quad \forall t \in
[0,T],
\end{equation}
i.e. $c_t(Q)$ is a $Q$-Potential.
\end{lemma}

\begin{proof}
By Theorem VII.8 of Dellacherie and Meyer \cite{DMeyer}, equality (\ref{c-A}) holds
true if $(c_t(Q))_{t \in [0,T]}$ is a positive $Q$-supermartingale
of class (D), that is $(c_{\sigma} (Q))_{\sigma \in S}$ is uniformly
integrable where $S$ is the family of all stopping times smaller or
equal to $T$.

The process $(c_t(Q))_{t \in [0,T]}$ is clearly adapted and positive
and, by hypothesis and from the cocycle property, $c_t(Q) \in L^1
(Q)$ for any $t \in [0,T]$. By the cocycle property we deduce that
for any $0 \leq s \leq t \leq T$
$$
E_Q[c_{t,T} (Q) \vert \mathcal{F}_s]=c_{s,T}(Q)-c_{s,t}(Q) \leq
c_{s,T}(Q),
$$
i.e. $(c_t(Q))_{t \in [0,T]}$ is a $Q$-supermartingale. Furthermore,
$c_T (Q)=0$. It remains to show that $(c_t(Q))_{t \in [0,T]}$ is of
class (D). This proof is postponed to the Appendix.
\end{proof}

\begin{remark}
Since in our setting $(c_t(Q))_{t \in [0,T]}$ is c\`{a}dl\`{a}g (see the Appendix for the proof), as
a particular case of the previous Lemma it follows that equation
(\ref{c-A}) holds for a c\`{a}dl\`{a}g $(A_t)_{t \in [0,T]}$.
\end{remark}

\medskip
Note that from (\ref{c-A}) it follows that $c_{t,u} (Q)=E_Q \left[
A_u - A_t \vert \mathcal{F}_t \right]$ for any $0 \leq t \leq u \leq
T$. Furthermore, the assumption $c_t(P)=0$ implies that for $Q=P$ we
have $A=A^P=0$.

\begin{lemma} \label{lemma8}
Let $\sigma, \tau$ be two stopping times such that $0 \leq \sigma
\leq \tau \leq T$ and $Q^1, Q^2$ be two probability measures
equivalent to $P$. Denote by $A^1,A^2$ the corresponding increasing
processes as in (\ref{c-A}).

Let $Q$ be the probability measure induced by
$$
q= \left\{
\begin{array}{rl}
q^1, \quad &\mbox{ on } H^1 = \sintr 0, \sigma  \sintr \cup
\sintr \tau, T \sintr \\
q^2, \quad &\mbox{ on } H^2 = \sintr \sigma, \tau  \sintr
\end{array}
\right.
$$
and denote by $A$ the corresponding process as in (\ref{c-A}).

Then
\begin{equation} \label{AH}
dA=dA^1 \vert _{H^1} +dA^2 \vert _{H^2} =1_{H^1} dA^1 +1_{H^2} dA^2.
\end{equation}
\end{lemma}

\begin{proof} Consider an arbitrary $t \in [0,T]$.

For $t \geq \tau$: we have that $c_t (Q)=c_t ^1 (Q^1)=E_{Q^1}
\left[A_T ^1 -A_t ^1 \vert \mathcal{F}_t \right]$.

For $\sigma \leq t < \tau$: from the cocycle property  we deduce
that
$$
\begin{array}{ll}
c_t (Q)&=c_{t,\tau} (Q)+E_Q \left[c_{\tau, T} (Q) \vert
\mathcal{F}_t \right] \\
&=E_{Q^2} \left[A_{\tau} ^2 -A_t ^2 \vert \mathcal{F}_t
\right]+E_{Q^2} \left[E_{Q^1}[A_{T} ^1 -A_{\tau} ^1 \vert
\mathcal{F}_{\tau} ] \vert \mathcal{F}_t \right] \\
&=E_{Q} \left[A_{\tau} ^2 -A_t ^2 +A_T ^1 -A_{\tau} ^1 \vert
\mathcal{F}_t \right]\\
&=E_{Q} \left[ \int _{(t,T]} \left(1_{H^1} dA^1 +1_{H^2} dA^2
\right) | \mathcal{F}_t \right].
\end{array}
$$

For $t \leq \sigma$: from the cocycle property and from the case
above we deduce that
$$
\begin{array}{ll}
c_t (Q)&=c_{t,\sigma} (Q)+E_Q \left[c_{\sigma, T} (Q) \vert
\mathcal{F}_t \right] \\
&=E_{Q^1} \left[A_{\sigma} ^1 -A_t ^1 \vert \mathcal{F}_t \right]
\\
&\quad +E_{Q^1} \left[ E_{Q} \left[\int _{\sintr \sigma,T
\sintr} \left(1_{H^1} dA^1 +1_{H^2} dA^2 \right)
| \mathcal{F}_{\sigma} \right] \vert \mathcal{F}_t \right]\\
&=E_{Q} \left[\int _{(t,T]} \left(1_{H^1} dA^1 +1_{H^2} dA^2
\right) | \mathcal{F}_t \right].
\end{array}
$$

Since $A_t \triangleq \int _{(0,t]} \left(1_{H^1} dA^1 +1_{H^2} dA^2
\right)$ is c\`{a}dl\`{a}g, predictable and increasing, $(A_t)_{t
\in [0,T]}$ is the process associated to $Q$ in the sense of
(\ref{c-A}). \end{proof}

\begin{corollary} \label{corollary9}
Let $\sigma _1, \sigma_2, \ldots , \sigma_n, \tau _1, \tau _2,
\ldots, \tau _n$ be stopping times such that $0 \leq \sigma_1 \leq
\tau _1 \leq \sigma _2 \leq \tau _2 \leq \ldots \leq \sigma _n \leq
\tau _n \leq T$ and let $Q$ be a probability measure equivalent to
$P$ and whose corresponding increasing process is denoted by $A$.
Set
\begin{equation} \label{H}
H \triangleq \sintr \sigma_1, \tau _1  \sintr \cup \sintr \sigma _2,
\tau _2 \sintr \cup \ldots \cup \sintr \sigma_n, \tau _n \sintr.
\end{equation}

Let $Q^H$ be the probability measure induced by $ q^H= q 1_H$ and
denote by $A^H$ the corresponding process as in (\ref{c-A}). Then
\begin{equation} \label{AHP}
dA^H=1_{H} dA.
\end{equation}
\end{corollary}

\begin{proof} The proof of this result is a repeated application of Lemma
\ref{lemma8} (with $Q^1=P$ and $Q^2 =Q$).
\end{proof}

\begin{lemma} \label{lemma10}
Let $Q$ be a probability measure equivalent to $P$ and $A$ be the
associated increasing process.

Then there exists a sequence $(\tau ^n)_{n \in \Bbb N}$ of stopping
times such that

(i) $\frac{dQ^{\sintr 0, \tau ^n \sintr}}{dP} \rightarrow_n ^{L^1}
\frac{dQ}{dP}$, where $Q^{\sintr 0, \tau ^n \sintr}$ denotes the
probability measure induced by $q^{\sintr 0, \tau ^n \sintr}=q
1_{\sintr 0, \tau ^n \sintr} $;

(ii) $c_{0,T} (Q^{\sintr 0, \tau ^n \sintr}) \uparrow c_{0,T} (Q)$;

(iii) $A_{\tau ^n}$ is bounded.
\end{lemma}

\begin{proof}
For any $n \in \Bbb N$ set $\sigma ^n \triangleq \inf \{ t\geq 0:
A_t \geq n \}$. Hence $\sigma ^n$ is a predictable stopping time.
For  any fixed $n$, take now a sequence $(\tau ^{n,m})_{m \in \Bbb
N}$ such that $\tau ^{n,m}$ is increasing (in $m$), $\tau ^{n,m} <
\sigma ^n$ on $\{ \sigma ^n >0 \}$ and $\tau ^{n,m} \uparrow \sigma
^n$. By definition of $\sigma ^n$ and from $\tau ^{n,m} < \sigma ^n$
it follows that $A_{\tau^{n,m}} \leq n$.

For any $\varepsilon>0$ small enough, take now $n$ and consequently
$m$ big enough to have $\Vert \frac{dQ^{\sintr 0, \tau ^{n,m}
\sintr}}{dP}-\frac{dQ}{dP}\Vert _1 \leq \varepsilon$. For such
indexes set $\tau ^{(n)}\triangleq \tau ^{n,m}$.

Take now $\tau ^n \triangleq \max _{k \leq n} \tau ^{(k)}$. It can
be checked that $(\tau ^n)_{n \in \Bbb N}$ is an increasing sequence
of stopping times and that $A_{\tau ^n} \leq n$ (since also $\tau ^n
< \sigma ^n$). Furthermore, since $\sigma ^n=T$ for sufficiently big
$n$ and $\tau ^n \uparrow T$, property (i) follows. Property (ii)
can be checked as usual (see, for instance, the proof of Lemma
\ref{lemma6}).
\end{proof}

\smallskip
\begin{lemma} \label{lemma11} Let $Q$ be a probability measure equivalent to $P$ and let $A$ be the associated
increasing process. Suppose that $A$ is bounded. Let $H$ be a
predictable set.

Suppose that $\mathcal{E} ( q 1_H \cdot B)$ is a uniformly
integrable martingale. Set $\frac{dQ^H}{dP}\triangleq \mathcal{E} (
q 1_H \cdot B)_T$ and denote by $A^H$ the associated increasing
process.

Then
\begin{equation} \label{ineqAH}
dA^H \leq dA,
\end{equation}
hence $A_T ^H \leq A_T$.
\end{lemma}

\begin{proof} First of all, we recall that the sets of the same form as
in (\ref{H}) form an algebra $\mathcal{A}$ (Boolean algebra) and
that the $\sigma$-algebra $\mathcal{P}$ of predictable sets is
generated by $\mathcal{A}$.

Consider now any predictable set $H \in \mathcal{P}$ satisfying the
hypothesis above. If $H \in \mathcal{A}$, we already know that $dA^H
= 1_H dA$ (from Corollary \ref{corollary9}). We will consider now
the general case.

Consider two stopping times $\sigma, \tau$ such that $0 \leq \sigma
\leq \tau \leq T$ and take a sequence $(H^n)_{n \in \Bbb N}
\subseteq \mathcal{A}$ such that
\begin{equation} \label{Hn}
\begin{array}{ll}
&E_Q \left[\int _0 ^T \vert 1_{H^n} -1_H \vert dA \right]\rightarrow
_n 0 \\
&E \left[\int _0 ^T \vert 1_{H^n} -1_H \vert dt \right]\rightarrow
_n 0.
\end{array}
\end{equation}
Set $Q^{H^n}$ the probability measure induced by $q^n = q 1_{H^n}$
and denote by $A^{H^n}$ the associated increasing process. Again
from Corollary \ref{corollary9} it follows that $dA^{H^n} =1_{H^n}
dA$ (since $H^n \in \mathcal{A}$). By (\ref{Hn}) we have that
$\frac{dQ^{H^n}}{dP} \rightarrow_n ^{L^1} \frac{dQ^H}{dP}$.

By lower semi-continuity of $c$ and by (\ref{c-A}), we get
$$
\begin{array}{rl}
E_{Q^H} \left[A_ {\tau} ^H - A _{\sigma} ^H \vert \mathcal{F}_
{\sigma} \right]&=c_{\sigma, \tau} (Q^H) \\
&\leq \liminf _n c_{\sigma, \tau} (Q^{H^n}) \\
&= \liminf _n E_{Q^{H^n}} \left[A_ {\tau} ^{H^n} - A _{\sigma}
^{H^n} \vert \mathcal{F}_ {\sigma} \right].
\end{array}
$$
Since $\int _{\sintr \sigma, \tau \sintr} 1_{H^n} dA \rightarrow _n
\int _{\sintr \sigma, \tau \sintr} 1_H dA$, $\int _{\sintr \sigma,
\tau \sintr} 1_{H^n} dA$ is uniformly bounded and $\mathcal{E}
(1_{\sintr \sigma, \tau \sintr \cap H^n} q \cdot B) \rightarrow _n
^{L^1} \mathcal{E} (1_{\sintr \sigma, \tau \sintr \cap H} q \cdot
B)$, then
\begin{equation} \label{Hn2}
E_{Q^{H^n}} \left[\int _{\sintr \sigma, \tau \sintr} 1_{H^n} dA
| \mathcal{F}_ {\sigma} \right] \rightarrow _n E_{Q^H}
\left[\int _{\sintr \sigma, \tau \sintr} 1_{H} dA |
\mathcal{F}_ {\sigma} \right].
\end{equation}

From (\ref{Hn}) and (\ref{Hn2}) it follows that
$$
\begin{array}{ll}
E_{Q^H} \left[A_{\tau} ^H - A_{\sigma } ^H \vert
\mathcal{F}_{\sigma} \right]&=E_{Q^H} \left[\int _{\sintr
\sigma,
\tau \sintr} dA^H | \mathcal{F}_ {\sigma} \right] \\
&\leq \liminf _n E_{Q^{H^n}} \left[A_{\tau} ^{H^n} - A_{\sigma }
^{H^n} \vert \mathcal{F}_{\sigma} \right] \\
&=E_{Q^H} \left[\int _{\sintr \sigma, \tau \sintr} 1_{H} dA
| \mathcal{F}_ {\sigma} \right],
\end{array}
$$
hence $E_{Q^H} \left[\int _{\sintr \sigma, \tau \sintr} dA^H
| \mathcal{F}_ {\sigma} \right] \leq E_{Q^H} \left[\int
_{\sintr \sigma, \tau \sintr} 1_{H} dA | \mathcal{F}_ {\sigma}
\right]$.

The same inequality holds if we replace $\sintr \sigma, \tau \sintr$
with any element $K \in \mathcal{A}$ (it is sufficient to sum over
intervals of the same form as in (\ref{H})), that is
\begin{equation} \label{H2}
E_{Q^H} \left[\int _0 ^T 1_K dA^H | \mathcal{F}_ {\sigma}
\right] \leq E_{Q^H} \left[\int _0 ^T 1_K 1_{H} dA |
\mathcal{F}_ {\sigma} \right].
\end{equation}
Moreover, by passing to the limit we obtain that inequality
(\ref{H2}) holds true for any $K \in \mathcal{P}$. So we get $dA^H
\leq 1_H dA$ as stochastic measures on $(0,T]$, hence $A_T^H \leq
A_T$. \end{proof}

\begin{lemma} \label{lemma12} Let $Q$ be a probability measure equivalent to $P$
and suppose that the corresponding increasing process $A$ is
bounded.

If $H^n$ is predictable, $H^n \uparrow (0,T] \times \Omega$ and
$Q^{H^n}$ is the probability measure induced by $q^{H^n}=q 1_{H^n}$,
then
$$
c_{0,T} (Q^{H^n}) \rightarrow _n c_{0,T} (Q) .
$$

\end{lemma}

\begin{proof} We already know (by Lemma \ref{lemma11}) that
\begin{equation} \label{eqA}
dA^{H^n}\leq 1_{H^n} dA.
\end{equation}
From $\frac{dQ ^{H^n}}{dP} \rightarrow _n ^{L^1} \frac{dQ}{dP}$,
inequality (\ref{eqA}) and lower semi-continuity of $c_{0,T}$ we get
$$
\begin{array}{ll}
c_{0,T}(Q) &\leq \liminf _n c_{0,T} (Q^{H^n})\\
&=\liminf_n E_{Q^{H^n}} \left[A_T ^{H^n} \right]
\\
&\leq \liminf_n E_{Q^{H^n}} \left[\int _{(0,T]} 1_{H^n} dA \right]
\end{array}
$$
Since $\int _{(0,T]} 1_{H^n} dA$ is bounded and $\frac{dQ
^{H^n}}{dP} \rightarrow _n ^{L^1} \frac{dQ}{dP}$, we have that
$$
\begin{array}{ll}
c_{0,T}(Q) &\leq \liminf _n c_{0,T} (Q^{H^n}) \leq \liminf_n
E_{Q^{H^n}} \left[\int _{(0,T]} 1_{H^n} dA \right] \\
&=E_Q \left[\int _{(0,T]} dA  \right]=c_{0,T} (Q),
\end{array}
$$
hence $c_{0,T} (Q^{H^n}) \rightarrow _n c_{0,T} (Q)$. \end{proof}

\begin{theorem} \label{theorem14} Let $Q$ be a probability measure
equivalent to $P$ and let $A$ be the associated increasing process.
Then there exists a sequence $(Q^n)_{n \in \Bbb N}$ of probability
measures with $q^n$ bounded such that $\frac{dQ ^{n}}{dP}
\rightarrow _n ^{L^1} \frac{dQ}{dP}$ and $c_{0,T} (Q^{n})
\rightarrow _n c_{0,T} (Q)$.
\end{theorem}

\begin{proof} From the arguments above (and by stopping arguments) we may
suppose $A$ bounded.

For any $n \in \Bbb N$ take $H^n \triangleq \{ \Vert q \Vert \leq
n\}$ and set $Q^n$ the probability measure induced by $q^n = q
1_{H^n}$. Hence $H^n$ is predictable and $H^n \uparrow (0,T] \times
\Omega$, it satisfies the hypothesis of Lemma \ref{lemma12}. It
follows that $\frac{dQ ^{n}}{dP} \rightarrow _n ^{L^1}
\frac{dQ}{dP}$ and (by Lemma \ref{lemma12}) that $c_{0,T} (Q^{n})
\rightarrow _n c_{0,T} (Q)$. \end{proof}

\bigskip

\bigskip
We are now ready to prove the representation of the penalty term $c$
in terms of $f$ (see Theorem \ref{thm3}).

\medskip
\begin{proof} \textbf{(of Theorem \ref{thm3})} Since (ii) is a straightforward consequence of (i) and of the
representation in (\ref{repres}), it remains to show that
\begin{equation} \label{eqC}
c_{0,T}(Q)=E_Q \left[ \int _0 ^T f(t, q_t)dt \right].
\end{equation}
By Lemma \ref{lemma4}, we already know that $c_{0,T} (Q) \leq E_Q
\left[ \int _0 ^T f(u,q_u) du \right]$ for any probability measure
$Q \sim P$.

Suppose that $\int _0 ^T f(t,\omega, q_t)dt \in L^1 (Q)$. For any $n
\in \Bbb N$ set $\sigma _n \triangleq \inf \{ t \geq 0 : \int _0 ^t
f(u,q_u) \geq n \}$. $(\sigma _n)_{n \geq 0}$ is a sequence of
stopping times such that $\sigma _n \uparrow T$.

Take now a sequence $(Q^m)_{m \in \Bbb N}$ of probability measures
as in Theorem \ref{theorem14}. Then
$$
\begin{array}{ll}
c_{0,T}(Q) &\leq E_Q \left[ \int _0 ^T f(u,q_u) du \right] \\
&= \lim _n E_Q \left[ \int _0 ^{\sigma _n} f(u,q_u) du \right] \\
&\leq \sup _n \lim _m E_{Q^m} \left[ \int _0 ^{\sigma _n}
f(u,q_u) 1_{\Vert q \Vert \leq m} du \right] \\
&\leq \lim_m \sup _n E_{Q^m} \left[ \int _0 ^{\sigma _n} f(u,q_u)
1_{\Vert q \Vert \leq m} du \right] \\
&=\lim _m E_{Q^m} \left[ \int _0 ^T
f(u,q_u) 1_{\Vert q \Vert \leq m} du \right] \\
&= \lim _m c_{0,T} ^m (Q^m)=\lim _ m c_{0,T} (Q^m)=c_{0,T}(Q),
\end{array}
$$
where the last equality is due to Theorem \ref{theorem14}. Equality
(\ref{eqC}) has therefore been established for $\int _0 ^T f(t,
q_t)dt \in L^1 (Q)$.

If $\int _0 ^T f(t,\omega, q_t)dt \notin L^1 (Q)$, by Fatou's Lemma
we get
$$
\begin{array}{ll}
c_{0,T} (Q) &\leq E_Q \left[ \int _0 ^T f(t, q_t)dt \right]
\\
&\leq \liminf_m E_{Q^m} \left[ \int _0 ^T f(t, q_t)1_{\Vert q
\Vert \leq m} dt \right] \\
& = \liminf _m c_{0,T} ^m (Q^m) \\
&=\liminf _m c_{0,T} (Q^m)= c_{0,T} (Q),
\end{array}
$$
hence $c_{0,T}(Q)=E_Q \left[ \int _0 ^T f(t, q_t)dt
\right]=+\infty$. The representation of $c_{s,t} (Q)$ (hence of
$c_{\sigma,\tau} (Q)$) can be deduced as usual. \end{proof}

\bigskip
\bigskip
\textbf{Acknowledgements} The authors thank two anonymous referees
for useful comments that improved this paper.

\bigskip

\section{Appendix}

Let $Q$ be a probability measure equivalent to $P$ and such that
$c_{0,T}(Q)<+\infty$.

In the following, we will prove that $(c_{t,T}(Q))_{t \in [0,T]}$ is
of class (D) and that it admits a c\`{a}dl\`{a}g modification.

The following corollary of Lemma \ref{lemma6} will be useful later.

\begin{corollary} \label{corollA1}
$\sup \{E_Q[c_{\tau,T}(Q)] \vert \tau \mbox{ stopping time s.t. } \mbox{ } P(\tau <T)\leq \frac{1}{n}\}\to _n 0$
\end{corollary}

The following result is a straightforward consequence of the cocycle property of $c$.

\begin{lemma} \label{lemma14} Denote by $S$ the family of all stopping times smaller or equal to $T$.

The family $(c_{\sigma,T}(Q))_{\sigma \in S}$ satisfies the following property: given any pair of stopping times $\sigma, \tau$ such that $0\leq \sigma \leq \tau\leq T$ then $c_{\sigma,T}(Q) \geq E_Q \left[c_{\tau,T}(Q) \vert \mathcal{F}_{\sigma} \right]$.
\end{lemma}

\begin{lemma} \label{lemma15}
The family $(c_{\sigma,T}(Q))_{\sigma \in S}$ is $Q$-uniformly
integrable.
\end{lemma}

\begin{proof}
We have to prove that
\begin{equation} \label{UI}
\lim _{n \to + \infty} \sup _{\sigma \in S} \int _{c_{\sigma, T}(Q) >n} c_{\sigma,T} (Q) dQ =0.
\end{equation}
Consider an arbitrary stopping time $\sigma \in S$ and set
$$
\sigma ^{(n)} =\left\{
\begin{array}{rl}
\sigma;& \quad \mbox{ if } c_{\sigma,T}(Q)>n \\
T;& \quad \mbox{ if } c_{\sigma,T}(Q)\leq n
\end{array}
\right.
$$
By the cocycle property we get
$$
\begin{array}{rl}
c_0 (Q) &=c_0 (Q^{\sigma ^{(n)}})+E_Q [c_{\sigma ^{(n)},T}(Q)] \\
&\geq E_Q [c_{\sigma ^{(n)},T}(Q)] \\
&=\int _{c_{\sigma,T}(Q)>n} c_{\sigma,T}(Q) dQ \geq n P(c_{\sigma,T}(Q)>n).
\end{array}
$$
Hence $P(c_{\sigma,T}(Q)>n) \leq \frac{c_0(Q)}{n}$ uniformly in $\sigma$, so we get
$$
\begin{array}{rl}
0&\leq \sup_{\sigma \in S} \int _{c_{\sigma,T}(Q)>n} c_{\sigma,T}(Q) dQ \\
&\leq \sup \{E_Q [c_{\tau,T} (Q)] \vert \tau \mbox{ stopping time s.t. } P(\tau <T) \leq \frac{c_0(Q)}{n} \}
\end{array}
$$
Since the last term tends to 0 as $n \to + \infty$ by Corollary \ref{corollA1}, (\ref{UI}) follows.
\end{proof}

\begin{lemma} \label{lemma16}
Let $\varepsilon >0$ such that $E_Q[ - \xi] >c_0(Q)-\varepsilon$
with $\xi \in \mathcal{A}_{0,T}$. Then for any pair of stopping
times $\sigma, \tau$ such that $0\leq \sigma \leq \tau \leq T$ it
holds that:
$$
E_Q[c_{\sigma,\tau}(Q)]\leq E_Q[u_{\sigma}(\xi)-u_{\tau}(\xi)]+\varepsilon.
$$
\end{lemma}

\begin{proof}
By translation invariance of $(u_{t,T})_{t \in [0,T]}$ it follows that $u_{\tau,T}(\xi-u_{\tau,T}(\xi))=0$, hence $\xi-u_{\tau}(\xi) \in \mathcal{A}_{\tau,T}$. Furthermore, by time-consistency and translation invariance of $u$ and by $\xi \in \mathcal{A}_{0,T}$ it follows that $u_{\tau}(\xi)-u_{\sigma}(\xi)\in \mathcal{A}_{\sigma,T}$ and that $u_{\sigma}(\xi)\in \mathcal{A}_{0,\sigma}$.

The cocycle property or, equivalently, the decomposition property $\mathcal{A}_{0,T}=\mathcal{A}_{0,\sigma}+\mathcal{A}_{\sigma,\tau}+\mathcal{A}_{\tau,T}$ implies that
$$
\begin{array}{rl}
c_0(Q)&=E_Q[c_{0,\sigma}(Q)]+E_Q[c_{\sigma,\tau}(Q)]+E_Q[c_{\tau,T}(Q)] \\
&\geq E_Q[-u_{\sigma} (\xi)]+E_Q[u_{\sigma} (\xi)-u_{\tau}(\xi)]+E_Q[u_{\tau} (\xi)-\xi] \\
&\geq E_Q [-\xi] \geq c_0(Q)-\varepsilon,
\end{array}
$$
where the first inequality follows from $c_{t,T}(Q)=\esssup_{\xi \in \mathcal{A}_{t,T}} E_Q \left[-\xi \vert \mathcal{F}_t \right]$.
By proceeding as above we get
$$
\begin{array}{rl}
c_0(Q)&\geq E_Q[-u_{\sigma} (\xi)]+E_Q[u_{\sigma} (\xi)-u_{\tau}(\xi)]+E_Q[u_{\tau} (\xi)-\xi] \\
&\geq c_0(Q)-\varepsilon \\
&\geq E_Q[-u_{\sigma} (\xi)]+E_Q[c_{\sigma,\tau}(Q)]+E_Q[u_{\tau} (\xi)-\xi] -\varepsilon,
\end{array}
$$
hence $E_Q[c_{\sigma,\tau}(Q)]\leq E_Q[u_{\sigma} (\xi)-u_{\tau}(\xi)]+\varepsilon$.
\end{proof}

\smallskip
\begin{lemma} \label{lemma18}
Let $\sigma \in S$. If $\{\sigma_n \}_{n \in \Bbb N}$ is a sequence of stopping times such that $\sigma _n \downarrow \sigma$, then $E_Q [c_{\sigma,\sigma _n}(Q)] \to _n 0$.
\end{lemma}

\begin{proof}
Suppose by contradiction that $E_Q [c_{\sigma,\sigma _n}(Q)]$ does not tend to 0 as $n \to + \infty$. Hence there exists $\varepsilon >0$ such that $E_Q [c_{\sigma,\sigma _n}(Q)]\geq \varepsilon >0$ for any $n \in \Bbb N$. Take now $\xi \in \mathcal{A}_{0,T}$ such that $E_Q[-\xi]\geq c_0(Q)-\frac{\varepsilon}{2}$. Hence, by Lemma \ref{lemma16},
$$
E_Q[u_{\sigma} (\xi)-u_{\sigma _n}(\xi)]\geq E_Q [c_{\sigma,\sigma _n}(Q)]-\frac{\varepsilon}{2} \geq \frac{\varepsilon}{2}
$$
for any $n \in \Bbb N$. This leads to a contradiction since $(u_{t,T})_{t \in [0,T]}$ admits a c\`{a}dl\`{a}g version with $u_{\sigma _n} (\xi) \to _n u_{\sigma}(\xi)$ in $L^1 (Q)$ (see Lemma 4 of Bion-Nadal \cite{BionN}).
\end{proof}

\smallskip
By the cocycle property it is easy to deduce the following result from the one above.

\begin{corollary} \label{corollA2}
Let $\sigma \in S$. If $\{\sigma_n \}_{n \in \Bbb N}$ is a sequence of stopping times such that $\sigma _n \downarrow \sigma$, then $E_Q [c_{\sigma _n,T}(Q)] \to _n E_Q [c_{\sigma,T}(Q)]$.
\end{corollary}

\begin{lemma} \label{lemma19}
$(c_{t,T}(Q))_{t \in [0,T]}$ admits a c\`{a}dl\`{a}g modification. Furthermore, if $(\overline{c}_{t})_{t \in [0,T]}$ denotes this modification, for any stopping time $\sigma \in S$ it holds $c_{\sigma,T}(Q)=\overline{c}_{\sigma}$ a.s.
\end{lemma}

\smallskip
We remark that this ends the proof of the statement in the beginning
of the appendix.

\smallskip
\begin{proof} We already know that $(c_{t,T}(Q))_{t \in [0,T]}$ is a positive $Q$-supermartingale
(see the proof of Lemma \ref{lemma7}) and that for any sequence $\{t_n \}_{n \in \Bbb N}$ in $[0,T]$
and such that $t_n \downarrow t$ it holds $E_Q [c_{t _n,T}(Q)] \to _n E_Q [c_{t,T}(Q)]$
(by Corollary \ref{corollA2}). By Theorem 4 at page 76 of Dellacherie and Meyer \cite{DMeyer}
it follows that $(c_{t,T}(Q))_{t \in [0,T]}$ admits a c\`{a}dl\`{a}g modification.
This implies that for any stopping time $\sigma \in S$ taking rational values it holds
$\overline{c}_{\sigma}=c_{\sigma,T}(Q)$ a.s.. For a general stopping time $\sigma \in S$
there exists a sequence $\{\sigma_n \}_{n \in \Bbb N}$ of finite stopping times taking
rational values and such that $\sigma_n \downarrow \sigma$. Hence
\begin{equation} \label{eq-cadlag}
\lim_{n \to + \infty} c_{\sigma _n, T}(Q)=\lim_{n \to + \infty} \overline{c}_{\sigma _n}=\overline{c}_{\sigma}\quad a.s.
\end{equation}
where the last equality follows from the fact that $(\overline{c}_{t})_{t \in [0,T]}$ is c\`{a}dl\`{a}g.

It remains to prove that $c_{\sigma,T}(Q)=\lim_{n \to + \infty}
c_{\sigma _n, T}(Q)$. This proof is quite standard and we include it
for completeness. By the cocycle property it follows that
$(c_{\sigma_n,T}(Q), \mathcal{F}_{\sigma_n})_{n \in \Bbb N}$ is a
positive reversed $Q$-supermartingale (see Neveu \cite{Neveu}). By
Proposition V-3-11 of Neveu \cite{Neveu}, $c_{\sigma_n,T}(Q)$
converges as $n \to + \infty$ to a positive
$\mathcal{F}_{\sigma}$-measurable random variable $\eta$ and
$E_Q\left[c_{\sigma_n,T}(Q) \vert \mathcal{F}_{\sigma}\right] \to _n
\eta$ a.s.. Since $E_Q\left[c_{\sigma_n,T}(Q) \vert
\mathcal{F}_{\sigma}\right]\leq c_{\sigma,T}(Q)$, we get $\eta \leq
c_{\sigma,T}(Q)$. Furthermore, by $Q$-uniform integrability of
$(c_{\sigma_n,T}(Q))_{n \in \Bbb N}$ (see Lemma \ref{lemma15}) we
get
$$
E_Q [c_{\sigma,T}(Q)]=\lim_n E_Q [c_{\sigma_n,T}(Q)]=E_Q[\eta]
$$
where the first equality is due to Corollary \ref{corollA2}. By the arguments above it follows
that $\eta=c_{\sigma_n}(Q)$ a.s., hence the thesis.
\end{proof}

\end{document}